\documentclass[12pt,a4paper,oneside]{article}
\usepackage{amsmath,amstext,amssymb,amsopn,amsthm}
\usepackage[latin1]{inputenc}

\theoremstyle{plain}
\newtheorem{thm}{Theorem}[section]

\newtheorem{lem}[thm]{Lemma}

\newtheorem{cor}[thm]{Corollary}

\theoremstyle{definition}

\newtheorem{exmp}{Example}[section]

\theoremstyle{remark}
\newtheorem*{rem*}{Remark}

\def\QED{\hfill\vrule height 1.5ex width 1.4ex depth -.1ex \vskip20pt}
\newcommand{\ind}{\mbox{\rm 1\hspace{-0.04in}I}}
\newcommand{\smallboullet}{{\scriptstyle \,^\bullet \,}}
\title{On the stability of call/put option's prices
 in incomplete models under statistical estimations.}
 \date{}
\author{Vostrikova L.  LAREMA, D\'epartement de Math\'ematiques, \\
Universit\'e d'Angers, France}

\begin{document}
\maketitle

\begin{abstract} In exponential semi-martingale setting for risky asset we estimate the 
difference of prices of  options when initial physical measure $P$ and
 corresponding martingale measure  $Q$ change to $\tilde{P}$ and $\tilde{Q}$
respectively. Then, we estimate $L_1$-distance of  option's prices for corresponding parametric models with 
known and estimated parameters. The results are applied to exponential Levy models with 
special choice of martingale measure as Esscher measure, minimal entropy measure and
 $f^q$-minimal martingale measure. We illustrate our results by
 considering GMY and CGMY models.
\\

\noindent\it 2000 MS Classification \rm : Primary 60G07, 60G51, 62F10 ; Secondary 91 B70  
   \\
\noindent\it Key words and phrases\rm :  stability, call/put option price, incomplete model,
 minimal entropy martingale measure, Esscher measure, $f^q$-martingale measure.
\footnotetext{
This work was supported by DNIPRO grant  14198ZL.
  \\}

\end{abstract}

\section{Introduction}
We consider the following semi-martingale model of risky asset
$S=(S_t)_{t\geq 0}$:
$$S_t= S_0 \exp(X_t)$$
where $X=(X_t)_{t\geq 0}$ is a semi-martingale. Usually the law of
this semi-martingale depend on unknown parameter, say $\theta \in
\Theta $, where $\Theta $ is some space.  For exemple,
in Black-Scholes model we have:
$$X_t= (\mu-\sigma^2/2)t +\sigma W_t$$
where $W=(W_t)_ {t \geq 0}$ is a standard Wiener process, the parameter $\theta=
(\mu , \sigma )$ and $\Theta = \mathbb R \times \mathbb R ^ {+, \star}$.

 In Geometric Variance Gamma model (cf. \cite{GMY1},\cite{GMY2}), as well known, 
\begin{equation}
\label{1}
X_t= \mu \tau _t +\sigma W_{\tau _t}
\end{equation}
where $\mu \in \mathbb R$ , $\sigma >0$, $W=(W_t)_{t\geq 0}$ is again Wiener process
and $(\tau _t)_{t\geq 0}$ is, independent from $W$, Gamma process with
parameters $(1, \nu )$, $\nu >0$.  In this case $\theta = (\mu ,\sigma, \nu )$
and $\Theta =\mathbb R\times\mathbb R^{+,*}\times \mathbb R^{+,*}$.

In GMY model, as well known ( cf.\cite{GMY1}, \cite{GMY2})
the process $X$ has the same structure as in \eqref{1} but with 
$(\tau _t)_{t\geq 0}$ being Levy process with Levy measure
\begin{equation}
\label{2}
\nu (dx)=
\frac{C\exp (-Nx)\ind _{ \{ x>0\}}}{x^{1+\alpha }}dx
\end{equation}
where $ \alpha <2$, $C>0$ and $N\geq 0$. Then, obviously,  $\theta=(C,N, \alpha )$ and $\Theta =
\mathbb R^{+,*}\times\mathbb R^{+}\times]-\infty,2[$.

In CGMY model(cf.\cite{GMY2},\cite{CGMY2}) the process $X$ is simply  a Levy
process with the Levy measure

\begin{equation}
\label{3}
\nu (dx)=
\displaystyle\frac{C\exp (-Nx)\ind _{\{x>0\}}+ C\exp (-Mx)\ind _{\{x<0\}}}{x^{1+\alpha }}dx
\end{equation}\\
where $C,M,N$ are positive constants and $\alpha <2$. Then we
have $\theta= (C,M,N, \alpha)$.

We will also mention Hyperbolic Levy process
$X^{\theta}=(X^{\theta}_t)_{t\geq 0} $ which is often used in
modelisation because of its flexibility to fit the form of
one-dimensional distributions of log of returns (cf.\cite{E},\cite{SH}).
As well known, there exist several parametrisations of Hyperbolic
Levy processes. Under one of them, say  $\theta = (\alpha, \beta,
\delta,\mu)$, the one dimensionnal densities of
$X^{\theta}_1$ with respect to Lebesgue measure are given by:\\
\begin{equation}
f(x)= \frac{\sqrt{\alpha ^2 -\beta ^2}}{2\alpha\delta K_1(\delta
  \sqrt{\alpha ^2 -\beta ^2})}\exp (-\alpha\sqrt{\delta^2 +(x-\mu)^2}
+ \beta(x-\mu))
\end{equation}\\
where $\alpha>0$, $0\leq |\beta |<\alpha$, $\delta>0$, $\mu \in \mathbb R$ and
$K_1(\cdot)$ is a Bessel function of the third type of index 1.
We know (cf.\cite{SH}) that Levy measure of this
process is equal to:\\
\begin{equation}
\nu(dx)= \frac{\exp(\beta x)}{|x|}\left(\exp ^{-\alpha |x|} + \int
^{+\infty}_0 \frac{\exp(-\sqrt{2y+\alpha ^2} |x|)}{\pi^2 y
  (J_1^2(\delta \sqrt{2y}) + Y_1^2(\delta \sqrt{2y}))}dy \right)
\end{equation}\\
where $J_1(\cdot)$ and $Y_1(\cdot )$ are Bessel functions of the first
and second type of index 1.

The classical procedure of calculus of   call/put option price
$\mathbb C_T$ of maturity time $T$  consists to take payoff function
given by a continuous in the space $D([0,T])$ functional $g(\cdot)$,
then to select in
the set of equivalent martingale measures $\mathcal M(P)$, supposed
non-empty, a  "good" one, say $Q$, and to put:
$$\mathbb C_T=\mathbb E _Q(g(S)).$$ 

As we know, there exist  many approaches to choose a "good" martingale measure: it can be done using
the minimisation of the risk in $L^2$-sense( see \cite{FS},\cite{S}),
 using the minimisation of
Hellinger integrals (see \cite{CSL}, \cite{G}), it can be based on the minimisation of entropy
(see \cite{MI},\cite{MIF}, \cite{CS}), one can take minimal $f^q$-martingale measures (see
\cite{JKM}) or use Esscher measures (see \cite{K},\cite{MI} ) e.t.c.\\
We remark that since the law of $X^\theta$ depends on $\theta$, the
price $\mathbb C_T$ does it as well. To ajust the "good" value of
$\theta$ one perform then so called calibration which is equivalent, from
statistical point of view, to find a minimal distance
estimator or contrast estimator with very special contrast.
About the properties of these estimators see for instance \cite{B},
\cite{M},\cite{KUT},
\cite{SV} and references there. One can use also another approach and consider maximum likelihood
estimators or Bayesian estimators for the unknown parameters. The properties of these estimators were studied, for example, in
\cite{IH}, \cite{KUT}, the
conditions for weak convergence of these processes in terms of Hellinger processes can be found in
\cite {V2}, \cite{V3}, \cite{JS}.
When the density of the law of $X$ with respect to some majorating measure can not be expressed explicitely or when it is too
complicated, one can use moment estimators ( see \cite{IH}). In practice often  the combination of some statistical estimations
and some calibration procedure also is used.

Let $\hat{\theta}$ be an estimator of unknown parameter $\theta$. Then,  we replace $\theta$ in formulas for  $\mathbb C_T(\theta)$
 by its estimator $\hat{\theta}$ and it becomes $\mathbb C_T(\hat{\theta})$.
 So, it is important from point of view of stability of the procesure
to measure the distance between estimated $\mathbb C_T (\hat{\theta})$ and
"true" price $\mathbb C_T(\theta)$. In this paper we are 
interested  to evaluate $L^1$ distance between these quantities,
namely $\mathbb E^{\theta}|\,\mathbb C _T(\hat{\theta})-\mathbb C _T(\theta)\,|$ where the expectation is
taken with respect to ''physical''  measure $P_{\theta}$. We remark that in the same manner one
can obtain the estimation of  $\mathbb E^{\theta}[d(\mathbb C _T(\hat{\theta}),\mathbb C
_T(\theta))]$ with different possible choise of the distance $d$.
We  notice the importance of use of consistent estimators of
$\theta$ in this procesure. In fact, usually $\mathbb C _T(\theta)\neq\mathbb C _T(\theta^{'})$ for
$\theta\neq \theta^{'}$. If the sequence of estimators is not consintent,
 then under some mild conditions one can extract a subsequence
 $(\hat{\theta}^n)$ converging $P-a.s.$
to $\theta + \delta$ with $\delta \neq 0$. Then
  $\mathbb E^{\theta}|\,\mathbb C _T(\hat{\theta}^n)-\mathbb C _T(\theta)\,|$ will
converge to  $|\,\mathbb C _T({\theta}+\delta)-\mathbb C _T(\theta)\,|$ 
which is different from zero. It means that without arbitrage for initial model we can have asymptotic arbitrage
consequences due to estimation procedure if  $\mathbb C_T(\theta+\delta)\neq \mathbb C_T(\theta)$.

In this paper we consider only payoff functions $g$ verifying the condition  \eqref{pay}. But similar results can be obtained in more general cases. The paper is organized in the following way. In  §2 we give the results for binary model, i.e. for the parametric models with 
two values of parameter. The main result is presented in Theorem 1. In  Corollary 1 the case of the processes
with independent increments is considered. Then, in section 3 we give the results for general parametric model. The main results are
presented in  Theorem 2 and Corollaries 2,3.
Finally, we apply the results for Levy processes, and we consider different possibilities to choose a martingale measure, namely
as Esscher measure,  Minimal entropy martingale measure and $f^q$- minimal martingale measure.
It is shown that under conditions of Theorem 1 we obtain the estimation of the type  \eqref{14}.
Then, the results are applied to Geometric Variance Gamma  and CGMY models. 

\section{Results for binary statistical model}

We suppose that we are given with  a filtered canonical space of
cadlag functions $(\Omega, \mathcal F, \mathbb F)$ where $\mathbb F=(\mathcal F _t)_{t\geq 0}$ is the
right-continuous filtration such that $\mathcal F= \bigvee
_{t\geq 0} \mathcal F_t$ and $\mathcal F_0=\{\emptyset ,
\Omega\}$. Let $P$ and $\tilde{P}$ be two equivalent probability measures
on $(\Omega , \mathcal F)$ and we denote by $P_t$ and $\tilde{P}_t$
the restrictions of these measures on the $\sigma$-algebra $\mathcal
F_t$, $t\geq 0.$  In this setting the measures $P$ and $\tilde{P}$ correspond to the
laws of our semimartingale $X=(X_t)_{t\geq 0}$ under two fixed values
of parameter.  We suppose that $X$ has predictable representation
property with respect to $P$  and the caracteristics of
$X$ are $(B,C,\nu)$ and  $(\tilde{B},\tilde{C},\tilde{\nu })$
respectively.  We remark that since the measures $P$ and $\tilde{P}$ are equivalent,
$C=\tilde{C}$ (P-a.s.) and we have the representation property with respect to $\tilde{P}$.  For more details about caracteristics see \cite{JS}.

We suppose that there is only two assets. For simplicity we assume that the interest rate $r$ of the bond $B=(B_t)_{t\geq 0}$ is equal to zero, 
i.e. $B_t=1$, and that the risky asset
$S=(S_t)_{t\geq 0}$ is given by:
$$ S_t=S_0\exp (X_t)$$
with $S_0=1$.
To avoid technical difficulties we suppose that
the processes $\ind_{\{x>1\}}\exp(x)\star \nu$ and
  $\ind_{\{x>1\}}\exp(x)\star \tilde{\nu}$
have bounded variation  on finite intervals. This supposition implies
that $S=S_0 \exp (X)$ is a special semimartingale under $P$ and $\tilde{P}$.

As usual we denote by $||P -\tilde{P} ||$ the variation distance
between the measures $P$ and $\tilde{P}$, i.e.
$$||P -\tilde{P}|| = 2\sup_{A\in\mathcal F} |P(A) - \tilde{P}(A)|$$
We recall that
$$||P -\tilde{P}|| = \mathbb E_P \big| 1- \frac{d \tilde{P}}{ d P}\big|.$$ 

 Let  $\mathcal M(P)$ and $\mathcal M(\tilde{P})$ be the sets
of equivalent martingale measures which are supposed to be non-empty.
Let $g$ be measurable functional in $D([0,T])$.
We choose, then,  using some procedure, two martingale measures: $Q$ and $\tilde{Q}$ to
calculate call/put option prices: $\mathbb C_T$ and $\tilde{\mathbb C _T}$
of maturity time $T$:
$$\mathbb C_T= \mathbb E_Q [g(S)],\hspace{1cm}\tilde{\mathbb C _T}=\mathbb E_{\tilde{Q}} [g(S)].$$

We introduce also dual measures $Q'$ and $\tilde{Q'}$ (cf. \cite{EPS})
by:
\begin{equation}
\label{dual}
 \frac{dQ'_T}{dQ_T}= S_T,
\hspace{1cm}\frac{d\tilde{Q'}_T}{d\tilde{Q}_T}= S_T.
\end{equation}

We notice that since $S$ is a martingale with respect to martingale
measure $Q$, $S' =1/S$ is also martingale but with respect to $Q'$. The same is
true for $\tilde{S}' =1/\tilde{S}$  with respect to $\tilde{Q'}$.

So, the measures involved in calculation can be represented by the following diagrams
containing initial measure, martingale measure and dual measure:
\begin{equation}\label{meg}
P \rightarrow  Q \rightarrow Q' \,\,\,\mbox{and}
\,\,\,\tilde{P } \rightarrow  \tilde{Q }\rightarrow \tilde{Q' }
\end{equation}
\begin{lem}Let $g$ be measurable functional in $D([0,T])$ verifying:
\begin{equation}\label{pay}
|g(x)| \leq c|x_T| + d
\end{equation}
 where $c, d$ are positif constants. Then for call/ put
option's price corresponding to $g$ we have:
$$| \mathbb C_T - \tilde{\mathbb C _T}|\leq c || Q'_T -\tilde{Q}'_T|| +  d|| Q_T -\tilde{Q}_T||$$
where $|| \cdot ||$ is a variation distance between the restriction of the corresponding measures
on $\sigma $-algebra $\mathcal F_T$.
\end{lem}

\it Proof. \rm

We have:
$$| \mathbb C_T - \tilde{\mathbb C _T} | = |\mathbb E_Q [g(S)]- \mathbb E_{\tilde{Q}}[g(S)]| \leq 
\mathbb E_Q\big[\big( c S_T + d\big)\big|1- \frac{d\tilde{Q_T}}{d Q_T}\big|\big]$$

But using  \eqref{dual} we obtain:
$$\mathbb E_Q\big(S_T\big|1- \frac{d\tilde{Q_T}}{d Q_T}\big|\big)= ||Q'_T - \tilde{Q}'_T||$$
 and by definition
$$\mathbb E_Q|1-\frac{d \tilde{Q}_T}{dQ_T}|= ||Q_T - \tilde{Q}_T||$$
\QED

It is known (see\cite{V1},\cite{VV},\cite{JS}) that the behaviour of
variation distance is closely related to the Hellinger distance and
Hellinger processes. Let $h(\frac{1}{2}, Q, \tilde{Q})= (h_t(\frac{1}{2}, Q, \tilde{Q}))_{t\geq 0}$
be the Hellinger process of order $1/2$ for the measures $Q$ and $\tilde{Q}$.

\begin{lem} We have the following estimation for the variation distance via Hellinger processes:
for $\epsilon >0$:
\begin{equation}\label{h1}
|| Q_T-\tilde{Q}_T|| \leq 4 \,[E_Q \,h_T(\mbox{$\frac{1}{2}$}, Q, \tilde{Q})]^{1/2}
\end{equation}
\begin{equation}\label{h2}
|| Q_T-\tilde{Q}_T|| \leq 3\sqrt{2\epsilon} + 2Q( \,h_T(\mbox{$\frac{1}{2}$}, Q, \tilde{Q})\geq \epsilon \, )
\end{equation}

\end{lem}
\it Proof \rm
See \cite{JS} p. 279.
\QED

To obtain the expressions for Hellinger processes we need the results on caracteristiques of the
process $X$ with respect to mentionned above measures.
First of all we remark that since the measure $Q$ 
is absolutely continuous with respect to $P$,
$X$ is a semi-martingale with respect to this measure and 
Girsanov theorem permit us to
find the caracteristics of  $X$
under $Q$ (see \cite{JS}, p. 159):
$$\left\{\begin{array}{l}
B^Q = B + \beta^Q \smallboullet C + l\cdot (Y^Q - 1)\star\nu\\
C^Q = C\\
\nu^Q = Y^ Q \cdot \nu
\end{array}\right.$$
where $l(\cdot)$ is a truncation function and $\beta^Q$ and $Y^Q$ are
 predictable functions verifying the
following integrability condition: for all $t\geq 0$ and P-a.s.
\begin{equation}\label{ic1}
 ((\beta^Q)^2 \smallboullet C)_t + (|l\cdot (Y^Q-1)|\star \nu)_t <\infty.
\end{equation}
Here and further $\smallboullet$ denotes a Lebesgue-Stielties integral and $\star$ 
means the integration with respect to a random measure ( for the details see \cite{JS}).
In  the mentionned above situation we say that $(\beta^{Q}, Y^{Q})$ are Girsanov parameters to pass from $P$ to $Q$.

The measures $Q^{'}$ and $\tilde{Q}^{'}$ are also absolutely continuous
with respect to $P$. In the following lemma we give predictable
caracteristics of $X$ with respect to the measures $ Q^{'}, \tilde{Q}$
and $\tilde{Q}^{'}$ via the caracteristics of the measure $P$.

\begin{lem} a) The predictable caracteristics of $X$ with respect to the
  measure $Q'$  via $P$ are given by:

$$\left\{\begin{array}{l}
B^{'} = B + (1+ \beta^Q)\smallboullet C + l\,(e^x Y^Q - 1)\star\nu \\
C^{'} = C\\
\nu ^{'} = e^{x} \,\,Y^Q \cdot \nu
\end{array}\right.$$
where $l(\cdot)$ is a truncation function and $(\beta^Q, Y^Q)$ are Girsanov parameters to
pass from $P$ to $Q$.

b) The predictable caracteristics of $X$ with respect to the
measure $\tilde{Q}$  via $P$ are given by:
$$\left\{\begin{array}{l}
B^{\tilde{Q}} = B + (\beta+ \beta^{\tilde{Q}})\smallboullet C + l\,
(Y^{\tilde{Q}} Y - 1)\star\nu\\
C^{\tilde{Q}} = C\\
\nu^{\tilde{Q}} = Y^{\tilde{Q}}\,\,Y \cdot \nu
\end{array}\right.$$
where $(\beta^{\tilde{Q}},Y^{\tilde{Q}})$ and
$( \beta , Y )$ are Girsanov parameters which permit us to pass from $\tilde{P}$ to $\tilde{Q}$ and from $P$ to $\tilde{P}$
respectively.

c) The predictable caracteristics of $X$ with respect to the measure
$\tilde{Q}'$ via $P$ are given by:

$$\left\{\begin{array}{l}
B^{\tilde{Q}'} = B + (1 + \beta + \beta^{\tilde{Q}}) \smallboullet C +
l\, (e^{x} Y^{\tilde{Q}} Y - 1)\star\nu \\
C^{\tilde{Q}'} = C\\
\nu ^{\tilde{Q}'} = e^{x}\,\, Y^{\tilde{Q}} \,\,Y \cdot\nu
\end{array}\right.$$
\end{lem}

\it Proof.  \rm To prove this Lemma we use  \eqref{meg}. We denote by $Z= (Z_t)_{t\geq 0},\,\,$
$\tilde{Z}= (\tilde{Z}_t)_{t\geq 0},\,\,$$Z^{'}= (Z^ {'}_t)_{t\geq 0},\,\,$$\tilde{Z}^ {'}= (\tilde{Z}^{'}_t)_{t\geq 0}$,
the processes such that for $t\geq 0$ and $P-a.s.$
$$Z_t=\frac{dQ_t}{dP_t},\,\,\,\tilde{Z}_t=\frac{d\tilde{Q}_t}{d\tilde{P}_t},\,\,\,
Z^{'}_t=\frac{dQ^{'}_t}{dP_t},\,\,\,\tilde{Z}^{'}_t=\frac{d\tilde{Q}^ {'}_t}{d\tilde{P}_t},$$
and $Q_t,\,\,\tilde{Q}_t,\,\,\,Q^{'}_t,\,\,\tilde{Q}_t^{'}$ stand for the restrictions
of the corresponding measures to the $\sigma$-algebra $\mathcal F_t$. To prove a) we note that
for all $t\geq 0$ we have:
$$Z_t^{'}= \frac{dQ^{'}_t}{dP_t}=\frac{dQ^{'}_t}{dQ_t}\frac{dQ_t}{dP_t}= e^{X_t}Z_t.$$
According to Girsanov theorem (see \cite{JS}, p. 160) the Girsanov parameters $(\beta^{Q^{'}}, \,Y^{Q^{'}})$
are given by: for $t\geq 0$
\begin{equation}\label{bqp} 
\beta^{Q^{'}}_t= \frac{1}{Z_{t-}^{'}}\frac{d\langle Z^{' c}, X^{c}\rangle_t}{dC_t}
\end{equation}
where $Z^ {' c}$ and $X^c$ denote continuous martingale part of the corresponding processes.
Using Ito formula for the function $f(x,y)= e^x y$ we find that 
$$Z_t^{' c}= \int_0^ t e^{X_{s-}}Z_{s-}dX^c_s + 
\int_0^ t e^{X_{s-}}dZ^ c_{s}.$$
Using the same formula as \eqref{bqp} for $\beta ^Q$ we obtain ($P$-a.s.) that  $\beta_t^{Q^{'}}=\beta_t^{Q}+1$.

Again according to Girsanov theorem
\begin{equation}
Y^{Q^{'}}= M^P_{\mu}\left(\frac{Z^{'}}{Z^{'}_{-}}\, |\,\tilde{ \mathcal P}\right)
\end{equation}
where $\tilde{\mathcal P}= \mathcal P\times \mathcal B(\mathbb R^*)$ is $\sigma$-algebra of predictable sets in 
$\tilde{\Omega}= \Omega \times [0,T]\times \mathbb R^{*}$ and for measurable non-negative functions $W(\omega, t,x)$
on $\tilde{\Omega}$
$$M^P_{\mu}(W)_T= \mathbb E_P[(W\star \mu)_T]$$
with $\mathbb E_P$ being the expectation with respect to $P$. Then
$$M^P_{\mu}\left(\frac{Z^{'}}{Z^{'}_{-}}\, |\,\tilde{ \mathcal P}\right)=
M^P_{\mu}\left(e^{\Delta X}\frac{Z}{Z_{-}}\, |\,\tilde{ \mathcal P}\right)$$
and, since the function $e^{\Delta X}$ is $\tilde{\mathcal P}$-measurable, we obtain that the right-hand side of the
previous equality is equal ($P$-a.s.) to:
$$
e^x \,M^P_{\mu}\left(\frac{Z}{Z_{-}} \,|\,\tilde{ \mathcal P}\right)$$
and we have a).

For b), c) we first write the caracteristics ot $X$ with respect to $\tilde{P}$ via $P$:
$$\left\{\begin{array}{l}
\tilde{B} = B + \beta \smallboullet C + l\cdot (Y - 1)\star\nu\\
\tilde{C} = C\\
\tilde{\nu} = Y \cdot \nu
\end{array}\right.$$
  Now we take the Girsanov parameters $(\beta ^ {\tilde{Q}}, Y^{\tilde{Q}})$ to pass from $\tilde{P}$ to $\tilde{Q}$:

$$\left\{\begin{array}{l}
B^{\tilde{Q}} = \tilde{B} + \beta^{\tilde{Q}} \smallboullet C + l\cdot (Y^{\tilde{Q}} - 1)\star\tilde{\nu}\\
C^{\tilde{Q}} = C\\
\nu^{\tilde{Q}} = Y^{\tilde{Q}} \cdot \tilde{\nu}
\end{array}\right.$$

Putting together these two decompositions we obtain b). Then, using the same procedure as in the proof of a), we obtain c).
\QED

Now we give the expressions for Hellinger processes.
To avoid technical difficulties we suppose that $X$ has no fixed points of discontinuity and that
for $\nu$ we have a desintegration formula. In fact, these suppositions are not too restrictives.
In fact, from one hand, atom's part can be also estimated , and, from another hand, a desintegtation formula 
with respect to some predictable increasing process always exists (see \cite{JS}, p. 77).
We introduce the following integrability condition:
\begin{equation}
\label{ic}
\int_{\mathbb R^*}|e^x-1|( d\nu^{Q}+d\nu^{\tilde{Q}}) < \infty.
\end{equation}

\begin{lem}
Let $X$ be a process without fixed points of discontinuity with respect to $P$.
We assume that  there exists a kernel $K(dx,t)$ such that we have a desintegration formula:
\begin{equation}\label{deco}
d\nu= K(dx,t)dC_t
\end{equation}
where $C$ is predictable variation of continuous martingale part of $X$ if it is not zero, 
and some increasing predictable process if not.
We suppose that  \eqref{ic} holds.
Then the Hellinger processes of order 1/2 of the measures $P$ and $\tilde {P}$, $Q$ and $\tilde{Q}$,
$Q'$ and $\tilde{Q}'$ are given respectively by:
$$h(\mbox{$\frac{1}{2}$}, P, \tilde{P})= \frac{1}{8}(\beta)^2 \smallboullet C + \frac{1}{2}\left(1-\sqrt{Y}\right)^2\star \nu,$$
$$h(\mbox{$\frac{1}{2}$}, Q, \tilde{Q})= \frac{1}{8}(\beta^Q - \beta^ {\tilde{Q}}-\beta)^ 2 \smallboullet C
+\frac{1}{2} \left(\sqrt{Y^ Q}- \sqrt{Y^ {\tilde{Q}}\cdot Y}\right)^2\star \nu,$$
$$h(\mbox{$\frac{1}{2}$}, Q', \tilde{Q}')= \frac{1}{8}(\beta^Q - \beta^ {\tilde{Q}}-\beta)^ 2 \smallboullet C
+ \frac{\exp(x)}{2}\left(\sqrt{Y^ Q}- \sqrt{Y^ {\tilde{Q}}\cdot Y}\right)^2\star \nu.$$

In addition we have ($P\times \lambda _C$ -a.s.)
$$\beta^Q - \beta^ {\tilde{Q}} - \beta =(\exp (x)-1)\left(Y^
{\tilde{Q}} \cdot Y - Y^Q \right)\star K(dx, \cdot)$$
where $\lambda _C$ is a positive measure with the distribition function $C$.
\end{lem}

\it Proof. \rm
 To obtain the expressions for the Hellinger processes we take in account that
 the compensator of $X$ has no atoms, we  use
the caracteristics  given in Lemma 2 and the formula in \cite{JS}, p. 221
(see also \cite{MS} for $X$ being the proceses with independent increments).
Since $e^{X}$ is a martingale with respect to the measures $Q$ and $\tilde{Q}$ we can write again 
using Ito formula and Girsanov theorem that (cf.\cite{JS}), p.556): $P$-a.s.
$$B^{Q} +\frac{1}{2} C^{Q}+ (e^x-1-l(x))\star \nu ^{Q} =0,$$
$$B^{\tilde{Q}} +\frac{1}{2} C^{\tilde{Q}}+ (e^x-1-l(x))\star \nu ^{\tilde{Q}} =0.$$
Then $P$-a.s.
$$B^{Q}- B^{\tilde{Q}} +(e^x-1-l(x))(Y^Q-Y^{\tilde{Q}}\cdot Y)\star \nu  =0.$$
 Taking in account  \eqref{deco}
and that $P$-a.s.
$$B^{Q}- B^{\tilde{Q}}= (\beta^{Q} -\beta^{\tilde{Q}}-\beta)\smallboullet C + l\cdot(Y^{Q} -1)\star \nu
- l \cdot (Y^{\tilde{Q}}\cdot Y -1)\star \nu $$
we obtain that
$$(\beta^{Q} -\beta^{\tilde{Q}}-\beta)\smallboullet C + 
(e^x - 1)\cdot(Y^{Q} - Y^{\tilde{Q}}\cdot Y)K(dx,\cdot)\smallboullet C=0$$
and that $P\times \lambda _C$ -a.s.
$$(\beta^{Q} -\beta^{\tilde{Q}}-\beta) + 
(e^x - 1)(Y^{Q} - Y^{\tilde{Q}} \cdot Y)K(dx,\cdot)=0.$$
\QED
 
Let us introduce the processes $\rho(Q,\tilde{Q})$ and $\rho(P,\tilde{P})$ 
which are closely related with the Hellinger processes, namely with their integral part
with respect to the compensator of the jump measure of $X$:
for all $t\geq 0$
\begin{equation}\label{r1}
\rho_t(Q, \tilde{Q})= \int_0^t\int_{\mathbb R^*}\left(\sqrt{Y^{\tilde{Q}}} - \sqrt{Y^Q}\right)^ 2 d\nu,
\end{equation}
\begin{equation}\label{r2}
\rho_t(P, \tilde{P})= \int_0^t\int_{\mathbb R^*}\left(1 - \sqrt{Y}\right)^ 2 d\nu.
\end{equation}

For a given non-negative constants $a, k$ we put
\begin{equation}\label{r}
A= 4a \sup_{0\leq t\leq T}\int_{\mathbb R^*} |e^x -1| e^{kx} K(dx, t)
\end{equation}
and we suppose that this quantity is finite $P$-a.s.
We introduce the functions 
$$p(x)= \frac{A |e^x -1|}{4} + 1,\hspace{1cm} q(x) = \frac{A |e^x -1|}{4} + e^x$$
We introduce also the processes $U = (U_t)_{t\geq 0}$ and $V=(V_t)_{t\geq 0}$ by:
\begin{equation}\label{U}
U_t = \int_0^t\int_{\mathbb R^*}p(x)d\rho_s(Q, \tilde{Q}) + \int_0^t\int_{\mathbb R^*}
ae^{kx}p(x) d\rho_s(P, \tilde{P})
\end{equation}
\begin{equation}\label{V}
V_t =\int_0^t\int_{\mathbb R^*}q(x)d\rho_s(Q, \tilde{Q}) + \int_0^t\int_{\mathbb R^*}
ae^{kx}q(x)\, d\rho_s(P, \tilde{P})
\end{equation}

\begin{lem} We suppose that $Y^{Q}$ and $Y^{\tilde{Q}}$ are bounded by $a e^{kx}$ where $a,k$ are non-negative constants
satisfying  $A <\infty$( $P$-a.s.), and that  \eqref{ic} holds.
Then we have:
$$h_T(\mbox{$\frac{1}{2}$}, Q, \tilde{Q})\leq U_T,$$
$$h_T(\mbox{$\frac{1}{2}$}, Q^{'}, \tilde{Q}^{'})\leq V_T$$
where the processes $U$ and $V$ are given by the formulas  \eqref{U},  \eqref{V}.
\end{lem}

\it Proof. \rm
We begin with the estimation of $h_T(\frac{1}{2},Q,\tilde{Q})$. Using Lemma 3 and Lemma 4 we write:
$$h_T(\mbox{$\frac{1}{2}$},Q,\tilde{Q})= 
\frac{1}{8}\int^ T_0\left(\int_{\mathbb R^*}(e^x-1)((Y^{Q} - Y^{\tilde{Q}}\cdot Y)K(dx,\cdot)\right)^2 dC_s$$
 $$+
 \frac{1}{2} \int_0^T\int_{\mathbb R^*}\left(\sqrt{Y^{\tilde{Q}}\cdot Y} - \sqrt{Y^Q}\right)^ 2 d\nu $$
For the first term on the right-hand side we have  by Schwartz inequality:
$$\int_0^T\left(\int_{\mathbb R^*}(e^x-1)((Y_s^{Q} - Y_s^{\tilde{Q}}\cdot Y)K(dx,\cdot)\right)^2dC_s \leq
 A \int^ T_0\int_{\mathbb R^*}|e^x-1|(\sqrt{Y^{Q}} - \sqrt{Y^{\tilde{Q}}\cdot Y})^2 d\nu$$
 where $A$ is given by  \eqref{r}.
This leads to the following inequality:
\begin{equation}\label{est}
h_T(\mbox{$\frac{1}{2}$}, Q, \tilde{Q})\leq \int_0^T\int_{\mathbb R^*}(\,\frac{A}{8}| e^x-1| + \frac{1}{2}\,)
\,\left(\sqrt{Y^{\tilde{Q}}\cdot Y} - \sqrt{Y^Q}\right)^ 2 d\nu,
\end{equation}
Now we remark that
$$(\sqrt{Y^{Q}} - \sqrt{Y^{\tilde{Q}}\cdot Y})^2 \leq 2(\sqrt{Y^{Q}} - \sqrt{Y^{\tilde{Q}}})^2 
+ 2Y^{\tilde{Q}}(1 - \sqrt{Y})^2 $$
and that $Y^{\tilde{Q}}$ is bounded by $a e^{k x}$. Then from the inequality  \eqref{est} we obtain the first result.
The second result can be obtained in similar way.
\QED

\begin{thm}Suppose that $X$ is a process without fixed points of discontinuity under $P$. We assume that  \eqref{ic},  \eqref{deco}
hold and that $Y^{Q}$ and $Y^{\tilde{Q}}$ are bounded by $a e^{kx}$ where $a,k$ are constants satisfying $A <\infty$ ( $P$-a.s.).
Then for payoff function satisfying  \eqref{pay} we have:

$$|\mathbb C_T - \tilde{\mathbb C_T}|\leq 4c \left[E_Q\, U_T\right]^{1/2} +
4d \left[E_{Q'}\, V_T\right]^{1/2},$$
Moreover, for $\epsilon >0$,

$$|\mathbb C_T - \tilde{\mathbb C_T}|\leq 3\sqrt{2\epsilon}(c+d) +
2c\, Q\left(\,U_T \geq \epsilon \,\right) +
2d\,Q'\left(\,V_T\geq \epsilon \,\right)$$
where the processes $U$ and $V$ given by the formulas  \eqref{U},  \eqref{V}
and $Q$, $Q'$ are martingale and dual martingale measure for $P$.
\end{thm}

\it Proof. \rm
We combine the Lemmas 2.2 and 2.5 to obtain the result.
\QED

Let us introduce the function
$$f(x)= \frac{A}{2}|e^ x-1| + \max (1, e^ x)$$
and the process $R=(R_t)_{t\geq 0}$ such that
\begin{equation}\label{R}
R_t = \int_0^t\int_{\mathbb R^*} f(x)d \rho _s(Q,\tilde{Q}) + \int_0^t\int_{\mathbb R^*}
ae^{kx} f(x) d\rho_s(P, \tilde{P})
\end{equation}

\begin{cor} Suppose that  $X$ is a process with independent increments under $P$ and $\tilde{P}$.  Assume that the conditions of Theorem 2.6 are satisfied. If in addition under the measures $ Q, \tilde{Q}$
 the process $X$ remains the process with independent increments then  for payoff function satisfying  \eqref{pay} we have:
$$|\mathbb C_T - \tilde{\mathbb C_T}| \leq 3\sqrt{2}(c+d) \,\sqrt{R_T}$$
\end{cor}
\it Proof. \rm
Use Theorem 1 and the fact that the processes $\rho(Q, \tilde{Q})$ and $\rho (P, \tilde{P})$ are deterministic.
\QED

\section{Results for general statistical model}

We suppose that $(\Omega, \mathcal F, \mathbb F)$ is filtered space endowed by the equivalent measures $P_{\theta}, \theta \in\Theta ,$ where $\theta$ is unknown parameter. We suppose that for each $\theta\in\Theta$, there exists a martingale measure $Q_{\theta}$.
 We denote as before  by $\mathbb C_T(\theta)$ the price  of risky asset obtained under physical measure $P_{\theta}$.
 Let $\hat{\theta}$ be an estimator of $\theta$ and let $\mathbb C_T(\hat{\theta})$ be the result of the replacement in $\mathbb C_T(\theta)$
 of the unknown parameter $\theta$ by its estimator.

We denote by $(\beta ^{\theta}, Y^{\theta})$ the Girsanov parameters to pass from $P_{\theta}$ to $Q_{\theta}$ and we
introduce the processes $U(\theta , \theta ')$ and $V(\theta , \theta ')$ by the formulas  \eqref{U},  \eqref{V}
with replacement $P, Q$ by $P_{\theta}, \,Q_{\theta}$, and  $\tilde{P}, \tilde{Q}$ by $P_{\theta '}, \,Q_{\theta '}$ respectively.
As before we assume that $S_0=B_0=1$ and $r=0$.

\begin{thm}
Suppose that the conditions of Theorem 2.6 are satisfied for each pair of measures $P_{\theta}$ and $P_{\theta '}$, $\theta\neq \theta'$,
$\theta, \theta'\in \Theta$. Then   for payoff function satisfying  \eqref{pay} we have:\\\\
$\mathbb E^ {\theta} |\mathbb C_T(\hat{\theta}) - \mathbb C_T(\theta )| \leq 2(c+d)\,P_{\theta}\left(|\hat{\theta} - \theta |>\epsilon\right)+\\$
 $$4c \sup _{|\theta -\theta^{'}|\leq \epsilon}\,\left[ E_{Q_{\theta}}\,U_T(\theta,\theta^{'})\right]^{1/2} +
 4d \sup _{|\theta -\theta^{'}|\leq \epsilon}\,\left[ E_{Q'_{\theta}}\,V_T(\theta, \theta^{'})\right]^{1/2}.$$
Moreover, for any $\epsilon>0$ we have:\\\\
$\mathbb E^ {\theta} |\mathbb C_T(\hat{\theta}) - \mathbb C_T(\theta )| \leq 
2(c+d)\,P_{\theta}\left(|\hat{\theta} - \theta |>\epsilon\right)+3\sqrt{2\epsilon}\,(c+d)\\ $
 $$+2c\, \sup _{|\theta -\theta^{'}|\leq \epsilon}\,Q_{\theta}\left(\,U_T(\theta, \theta^{'})\geq \epsilon\,\right)
 + 2d \,\sup _{|\theta -\theta^{'}|\leq \epsilon} \,Q^{'}_{\theta}\left(\,
 V_T(\theta, \theta^{'})\geq \epsilon\,\right)$$
 where $Q_{\theta}$ is the martingale measure of ''physical'' measure  $P_{\theta}$ and $Q^{'}_{\theta}$ is the respective dual measure.
\end{thm}

 \it Proof. \rm
 We remark that
$$\mathbb E^ {\theta} |\mathbb C_T(\hat{\theta}) - \mathbb C_T(\theta )| =\int_{\Omega} |\mathbb C_T(\theta ^{'}) - \mathbb C_T(\theta )|d P_{\hat{\theta}}(\theta ^{'})$$
and that for any $\epsilon >0$ the right-hand side can be majorated by:
$$2\sup_{\theta\in\Theta} \mathbb C_T(\theta)\mathbb \cdot P\left( |\hat{\theta}-\theta| >\epsilon\right) + \sup_{|\theta -\theta^{'}|\leq \epsilon}
 |\mathbb C_T(\theta) - \mathbb C_T(\theta ^{'} )|$$
 Due to (\ref{pay}) and martingale properties of $S$, we have  $\mathbb C_T(\theta) \leq  c+d$.
 Then we use the estimations of Theorem 2.6 to conclude.
\QED

Let also
\begin{equation}\label{R1}
R_T(\theta , \theta ') = \int_0^T\int_{\mathbb R^*} f_{\theta , \theta '}(x)d\rho _s(Q_{\theta}, Q_{\theta '}) + \int_0^T\int_{\mathbb R^ *}
a_{\theta ,  \theta '}e^{k_{\theta,  \theta ' }x}\,f_{\theta , \theta '}(x) d\rho_s(P, \tilde{P})
\end{equation}
where $f_{\theta , \theta '}$ , $A_{\theta , \theta '}$, $a_{\theta , \theta '}$ and $k_{\theta ,\theta '}$ are the fuction and the constants corresponding to $f$,  $A$, $a$ and $k$ of Theorem 2.6.

\begin{cor}Suppose that the process $X$ is a process with independent increments under $P_{\theta}, \theta \in \Theta$, as well as
under corresponding martingale measures $Q_{\theta}, \theta \in \Theta$. Suppose also that the conditions of Theorem 3.1 are satisfied.
Then  for payoff function satisfying  \eqref{pay} we have:
$$ \mathbb E^{\theta} |\mathbb C_T(\hat{\theta}) - \mathbb C_T(\theta )| \leq 2(c+d)\,P_{\theta}\left(|\hat{\theta} - \theta |>\epsilon\right)+
3\sqrt{2}(c+d)\,\left[\,\sup _{|\theta -\theta^{'}|\leq \epsilon } R_T(\theta , \theta ')\,\right]^{1/2}$$
\end{cor}

\begin{cor}Suppose that we have a sequence of processes with independent increments involving the physical measures
$(P^n_{\theta})_{n\geq 1}, \,\theta\in \Theta,$ the corresponding martingale measures $(Q^n_{\theta})_{n\geq 1},\,\theta\in\Theta,$
 and the respective sequence of the consistent estimators $(\hat{\theta}^n)_{ n \geq 1}.$ Suppose also that the conditions of Theorem 3.1 are satisfied. Let $R^n_T(\theta , \theta ')$ be defined by  \eqref{R1} with replacement
of $P_{\theta},Q_{\theta}$ and   $P_{\theta '},Q_{\theta '}$ by  $P^n_{\theta},Q^n_{\theta}$ and  $P^n_{\theta '},Q^n_{\theta '}$ respectively.

If uniformly in the neighbourhood of $\theta$
  as $n\rightarrow \infty$
$$R^n_T(\theta, \theta ')\rightarrow 0$$
 then  for payoff function satisfying  \eqref{pay} we have:
$$\mathbb E^n_{\theta}|\mathbb C_T(\hat{\theta}^ n) - \mathbb C_T(\theta )|\rightarrow 0$$
where  $E^n_{\theta}$ is a mathematical expectation with respect to $P_{\theta}^ n$.
\end{cor}

\section{Applications to Levy processes}

Suppose now that $X$ is Levy process with parameters $(b,c,\nu)$ under the measure $P$.
We emphasize that here $\nu$ is no more the
compensator of the measure of jumps of $X$ but a Levy measure, i.e. positive $\sigma$-finite measure on $\mathbb R$ such that
$$\int_{\mathbb R^ *} (x^2\wedge 1)d\nu <\infty.$$
We recall that the caracteristic function of $X_t$ for $t\geq 0$ and $\lambda \in\mathbb R$ is given by:
$$\phi_t(\lambda)= \exp (t \psi(\lambda ))$$
where $\psi(\lambda)$ is a caracteristic exponent of Levy process,
$$\psi(\lambda)= i b \lambda - \frac{1}{2}\lambda ^ 2 c + \int_{\mathbb R^ *}(\exp(i\lambda x)-1-i\lambda l(x))d\nu ,$$
and $l$ is the truncation function.
Let now $\tilde{P}$ be the measure coresponding to the parameters $(\tilde{b},\tilde{c},\tilde{\nu})$. According to 
 Corollary 2.7 of section 2 we have to find, for chosen martingale measures $Q$ and $\tilde{Q}$, the Girsanov parameters 
$(\beta^Q, Y^Q)$ and $(\beta^{\tilde{Q}}, Y^{\tilde{Q}})$ and write the expressions for the processes $\rho(Q,\tilde{Q})$
and $\rho(P,\tilde{P})$. We recall that as before $S_t= \exp(X_t)$. Let $r$ be positive constant, and,
let us suppose that the value process of the bond is deterministic and given by
$B_t=\exp (r t).$
\subsection{Esscher measures}
Esscher measures play very important role in actuary theory as well as in the option pricing theory
and they were studied in  \cite{K},\cite{MI}, \cite{MIF}. Let $$ D= \{\lambda\in\mathbb R\,|\, E_P e^{\lambda X_1}<\infty\}$$ where $E_P$ is the expectation with respect to
the physical  measure $P$.
Then for $\lambda\in D$ we define Esscher measure $P^{ES}$ of the parameter $\lambda$ and risk process $(X_t)_{t\geq 0}$
by : for $t\geq 0$
$$\frac{dP^{ES}_t}{dP_t} = \frac{e^{\lambda X_t}}{ E_P[ e^{\lambda X_t}]}$$
It is known that $(e^{-rt}S_t)_{t\geq 0}$ is a martingale under $Q=P^{ES}$ iff 
$$\psi(-i(1+\lambda)) - \psi(-i \lambda )= r$$
and the last equation is equivalent to:
\begin{equation}\label{l1}
b+ (\frac{1}{2}+\lambda) c + \int_{\mathbb R^ *}(( e^x\,-1)\,e^{\lambda x} - l(x))d\nu = r
\end{equation}
About existence and uniqueness of solution of  \eqref{l1}  see \cite{HS} and  \cite{KS}.\\

Suppose again that $X$ is Levy process with parameters $(b,c,\nu)$ under $P$, and that it has the parameters
 $(\tilde{b},c,\tilde{\nu })$ under $\tilde{P}$. Suppose that the solution of   \eqref{l1} exists as well as the solution 
of the same equation with the replacement $(b, c, \nu)$ by $(\tilde{b},c,\tilde{\nu })$ denoted  $\lambda ^ {*}$ and $\tilde{\lambda}^{*}$
respectively. Then $Q= P^{ES}(\lambda^{*})$ and $\tilde{Q} = P^{ES}(\tilde{\lambda}^{*})$.

Now we show that  the Girsanov parameters for $Q$ and $\tilde{Q}$ are: $\beta^{Q}= \lambda^{*}\,\, Y^{Q}= e^{\lambda^{*}\,x}$ and
$\beta^{\tilde{Q}}= \tilde{\lambda}^{*},\,\, Y^ Q= e^{\tilde{\lambda}^{*}\,x}$ respectively.
We write 
\begin{equation}\label{l3}
Z_t= \frac{dQ_t}{dP_t} = \frac{e^{\lambda^{*}X_t}}{\phi(-i\lambda^{*}t)}
\end{equation}
From the formula  \eqref{l3} we see that
$$\frac{Z_t}{Z_{t-}} = e^{\lambda^{*} \Delta X_t}$$
and according to Girsanov theorem
$$Y^{Q}= M^ P_{\mu}\left(e^{\lambda^{*}\Delta X}|\tilde{\mathcal P}\right)= e^{\lambda^{*}x}.$$
We use Ito formula to find $Z^c$ :
$$Z^ c_t= \int^ t_0 \frac{\lambda^{*}\exp (\lambda^{*}X_{s-})}{\phi(-i\lambda^{*}s)} d X_s^c$$
and, hence,
$$\beta_t^ Q= \frac{1}{Z_{t-}}\frac{d\langle Z^ c, X^ c\rangle _t}{ dC _t} = \lambda^{*} .$$
Now, we have to write the expression of $\rho _T(Q,\tilde{Q})$ and $\rho _T(P,\tilde{P})$:
$$\rho _T(Q,\tilde{Q})=T\int_{\mathbb R^ *}(\sqrt{e^{\lambda^{*}x}} - \sqrt{e^{\tilde{\lambda}^{*}x}})^2 d\nu$$
$$\rho _T(P,\tilde{P})=T\int_{\mathbb R^ *}(1 - \sqrt{Y})^2 d\nu$$
where $Y= d\tilde{\nu}/d\nu.$

In the case when $\lambda^{*}\leq 0$ and $\tilde{\lambda}^{*}\leq 0$ we can find  easily that the conditions of Theorem 2.6
are verified with $k=0$ and $a=1$. We remark that mean value theorem gives:
$$(\sqrt{e^{\lambda^{*}x}} - \sqrt{e^{\tilde{\lambda}^{*}x}})^2 \leq |x|^2(\lambda^{*} -\tilde{\lambda}^{*})^2$$
So,  for payoff function satisfying  \eqref{pay} we obtain the estimation:
\begin{equation}\label{14}
|\mathbb C_T - \tilde{\mathbb C}_T| \leq  T (\lambda^{*} -\tilde{\lambda}^{*})^2 \int_{\mathbb R^ *}f(x)x^2 d\nu + T
\int_{\mathbb R^ *} f(x) (\sqrt{d\nu} - \sqrt{d\tilde{\nu}})^2
\end{equation} 
where $f(x)= \frac{A}{2}|e^ x-1| + \max (1, e^ x)$ and  $A= 4aT \int_{\mathbb R^ *}|e^ x-1| d\nu$ .
 In the case when $\lambda^{*}$ and/or  $\tilde{\lambda}^{*}$ are not negatif we can obtain similar estimations .
\subsection{Minimal entropy measures}

Let $Q$ and $P$ be two equivalent probability measures then the relative entropy of $Q$ with respect to $P$ 
( or  Kulback-Leibler information in $Q$ with respect to $P$) is:
$$H(Q|P) = E_Q \left(\ln (\frac{dQ}{dP})\right) = E_P \left(\frac{dQ}{dP} \ln (\frac{dQ}{dP})\right)$$
We are interested in minimal entropy martingale measure, i.e. the measure $P^{ME}$ such that $(e^{-rt}S_t)_{t\geq 0}$
is a $P^{ME}$-martingale, and that for all $Q$ martingale measures 
$$H(P^{ME}|P) \leq H(Q|P)$$
 It turns out (cf.\cite{MIF}) that in the case of Levy processes $P^{ME}$ is nothing else as Esscher measure but for another
risque process $(\hat{X}_t)_{t\geq 0}$, namely for the process appearing in the representation:
$$S_t= S_0\, \mathcal E (\hat{X})_t$$
where $\mathcal E(\cdot)$ is Dolean's-Dade exponential,
$$\mathcal E(\hat{X})_t=
 \exp (\hat{X}_t - \frac{1}{2}\langle \hat{X}\rangle _t) \prod_{0\leq s\leq t}(1+\Delta\hat{X}_s)e^{-\Delta \hat{X}_s}$$

Writing Ito formula for $f(x)=e^ x$ we obtain that  $S_t= S_{t-} d\hat{X}_t$ with
$$\hat{X}_t= X_t +\frac{1}{2} \langle X^ c\rangle _t + \int _0^ t \int_{\mathbb R^*}(e^x - 1 - x)d\mu (x)$$
where $\mu$ is the measure of jumps of $X$.
This permits us to find the caracteristics of $\hat{X}$:
$$\left\{\begin{array}{l}
\hat{B} = B + \frac{1}{2}\langle X^c\rangle +  (e^x- 1 - x)\star\nu \\
\hat{C}= C\\
\hat{\nu} = (e^ x-1) \cdot \nu
\end{array}\right.$$
 We see that if $X$ is a Levy process verifying $\int_{\mathbb R^ *}|e^ x - 1| d\nu < \infty$
where $\nu$ is a Levy measure of $X$, then $\hat{X}$ is also Levy process and the parameters of $\hat{X}$
are:
$$\left\{\begin{array}{l}
\hat{b} = b + \frac{1}{2}c +  (e^x- 1 - x)\star\nu\\
\hat{c}= c\\
\hat{\nu} = (e^ x-1) \cdot \nu
\end{array}\right.$$
 Now let
  $ D= \{\lambda\in\mathbb R\,|\, E_P e^{\lambda \hat{X}_1}<\infty\}$
and let us introduce Esscher measure corresponding to the risque process $\hat{X}$ and $\lambda \in D$ : for $t\geq 0$
$$\frac{dP^{ME}_t}{dP_t} = \frac{e^{\lambda \hat{X}_t}}{ E_P [e^{\lambda \hat{X}_t}]}$$
 We remark that one can write easily the caracteristic function of $\hat{X}$ and the expression for caracteristic exponent:
$$\hat{\psi}(\lambda)= i\lambda( b+\frac{1}{2} c + (e^x - 1 - x)\star \nu)- \frac{1}{2}\lambda ^ 2c +
 \int_{\mathbb R^ *}(\exp(i\lambda x)-1-i\lambda l(x))(e^ x - 1)d\nu .$$
where  $l$ is the truncation function.

As it was mentionned before, this measure is a martingale measure for $(e^{-rt}S_t)_{ŧ\geq 0}$ iff
$$\hat{\psi}(-i(1+\lambda)) - \hat{\psi}(-i\lambda) = r$$
and the last equation is equivalent to:
\begin{equation}\label{m11}
b+ (\frac{1}{2}+\lambda) c + \int_{\mathbb R^ *}(( e^x\,-1)\,e^{\lambda (e^x -1)}  - l(x))\,d\nu = r
\end{equation}
About existence and uniqueness of solution of  \eqref{l1}  see \cite{HS} and  \cite{KS}.

Let us suppose that the solution $\lambda^{*}$ of the equation  \eqref{m11} exists as well the solution $\tilde{\lambda}^{*}$ 
of the similar equation with replacing $(b,c,\nu)$ by $(\tilde{b},c,\tilde{\nu})$.
We can show in the same way as before that Girsanov parameters of minimal entropy martingale measures are $(\lambda^{*}, e^{\lambda^{*}(e^x-1)})$ and $(\tilde{\lambda}^{*}, e^{\tilde{\lambda}^{*}(e^x-1)})$ respectively.
Then,  if  $\lambda^{*}$ and $\tilde{\lambda}^{*}$ are negatifs,  for payoff function satisfying  \eqref{pay} we have:

\begin{equation}\label{m12}
|\mathbb C_T - \tilde{\mathbb C}_T|\leq T (\lambda^{*} -\tilde{\lambda}^{*})^2 \int_{\mathbb R^ *}f(x) (e^x-1)^2 d\nu +
 T \int_{\mathbb R^ *} f(x) (\sqrt{d\nu} - \sqrt{d\tilde{\nu}})^2
\end{equation} 
where  $f(x)= \frac{A}{2}|e^ x-1| + \max (1, e^ x)$ and $A= 4aT \int_{\mathbb R^ *}|e^ x-1|d\nu$ .
 In the case when $\lambda^{*}$ and/or  $\tilde{\lambda}^{*}$ are not negatif we can obtain similar estimations .\\

\begin{exmp}1
 In Geometric Variance Gamma model the parameters $(b,c,\nu)$ are equal to $(0,0,\nu)$. The Levy measure 
of this model has the following form:
$$\nu(dx)= \frac{C(\ind_{\{x<0\}}e^{-M|x|} +\ind_{\{x>0\}}e^{-N x})}{|x|}dx$$
where $C>0$ and $M, N \geq 0$. 

We denote the left-hand side of  \eqref{m11} with given $\nu$ by  $\hat{f}$.  It is known (see \cite{MIN}) that if $0\leq N\leq 1$, or $N>1$ and $\hat{f}(0)\geq r$, then $\lambda ^{*}<0$.
 If $N>1$ and $\hat{f}(0)<r$, then $\lambda^{*}$ does not exist. So, 
 we have the estimation  \eqref{m12} when the solution of \eqref{m11} exists.
\end{exmp}

\begin{exmp}
 In Geometric CGMY model the parameters $(b,c,\nu)$ are equal to $(0,0,\nu)$. The Levy measure 
of this model has the following form:
$$\nu(dx)= \frac{C(\ind_{\{x<0\}}e^{-M|x|} +\ind_{\{x>0\}}e^{-N x})}{|x|^{1+\alpha}}dx$$
where $\alpha <2$, $C>0$ and $M, N \geq 0$. We recall that the case of $\alpha=0$ corresponds to Geometric Variance Gamma model and it was already considered.

We denote again the left hand side of  \eqref{m11} by $\hat{f}.$ It is known (cf.\cite{MIN}) that if $M=N=0$ and $0<\alpha <2$ then $X$ is symetric stable process and if, in addition $C>0$, then $\lambda^{*}<0$.
 If $0\leq N\leq 1$ or if $N>1$ and $\hat{f}(0)\geq r$ then again $\lambda^{*}<0$ . If $N>1$ and $\hat{f}(0)<r$
the equation  \eqref{m11} has no solution. So,  we have  the estimation  \eqref{m12}  when the solution of \eqref{m11} exists.
\end{exmp}
\subsection{$f^q$- martingale measures}
These measures take part of the measures minimising so called $f$-divergence between two probability measures. Let $Q$ and $P$ be
two  probability measures, $Q<<P$, and $f$ be a convex function with the values in $\mathbb R^{+,*}$. Then $f$-divergence (cf.\cite{CS}) 
of $Q$ given $P$,
denoted $f(Q\,|\,P)$ is given by:
$$f(Q\,|\,P)= E_P[\,f\left(\frac{dQ}{dP}\right)]$$
If $f(x) = x \ln x$ we obtain as $f(Q\,|\,P)$ the entropy or Kulback-Leibler information, if $f(x)=|1-x|$ we obtain the variation distance,
if $f(x)= (1-x)^ 2$ we obtain  variance squared distance, if $f(x)= (1-\sqrt{x})^2$ we obtain Hellinger distance.
 We remark also that the minimisation of variance squared distance is equivalent
 to minimise
$E_P[(\frac{dQ}{dP})^2]$, and that the minimising of Hellinger distance is equivalent to minimise $- E_P(\sqrt{\frac{dQ}{dP}})$.

In the papers \cite{CS},\cite{CSL}, \cite{JKM} the authors consider $f$-divergences with
$$f(x) = \left\{
\begin{array}{lll}
-x^q, & \mbox{if}\,\,& 0<q<1,\\
x^q, & \mbox{if}\,\,&q<0\,\, \mbox{or}\,\, q>1.
\end{array}\right.$$
It is not difficult to see that such $f$ is a convex function. It was shown that in the case of Levy processes the Girsanov parameters
$(\beta _q, Y_q)$ of the measure $P^{(q)}$ minimising $f$-divergence  given by the above expression, are deterministic. So, 
$X$ is also Levy process under $P^{(q)}$. 

It can be also shown that if $X$ is not monotone Levy process and if we allow as $P^{(q)}$ not only equivalent, but also absolute continuous measures, then  the Girsanov parameters $(\beta _q, Y_q)$  are  unique 
minimizers of the function
$$k(\beta, Y)= \frac{q(q-1)}{2}\beta^ 2 c + \int_{\mathbb R^*}(Y^q -1-q(Y-1))d\nu$$
under constraint  $$b+c\beta  + \int_{\mathbb R^ *}(xY(x)-l(x))d\nu=0$$
on the set  $$\mathcal A = \{ (\beta, Y)\,|\,\beta\in \mathbb R, Y \geq 0,\,
\int_{\mathbb R^ *}|\,xY(x)-l(x)\,|d\nu < \infty\}$$

Via an application of the Kuhn-Tucker theorem it can be  shown that 
$$Y_q(x)= \left\{
\begin{array}{lll}
\left(1+(q-1)\beta_ q (e^ x -1)\right)^{\frac{1}{q-1}}&\mbox{if}&1+(q-1)\beta_ q (e^ x -1)\geq 0,\\
0 & &\mbox{in opposite case,}
\end{array}\right.$$
where  $\beta _q$ is the first Girsanov parameter which can be find from the constraint. We remark that if in addition
$$\mbox{supp}(\nu) \subseteq \{x: 1+(q-1)\beta_ q (e^ x -1)>0\}$$
 then $P^{( q)}$ is equivalent to $P$.  We will suppose that the last  condition is satisfied.\\

Let  $(\beta_q, Y_q)$ and $(\tilde{\beta}_q, \tilde{Y}_q)$ be Girsanov parameters of $f^ q$-minimal martingale measures for $P$ and $\tilde{P}$ respectively.
 To evaluate $\rho(P^{(q)},\tilde{P}^{(q)})$ we remark that  
$$(\sqrt{Y_q(x)} - \sqrt{\tilde{Y}_q(x)})^ 2 \,\leq \, C \,(e^x-1)^2 (\beta_{q} - \tilde{\beta}_{q})^2$$
with some constant $C$.
So, we have the estimations similar to  \eqref{m12}.

\end{document}